\documentclass[12pt]{article}
\title{Symplectic Resolutions for Nilpotent Orbits (II)}
\author{Baohua Fu}
\usepackage{amsmath,amssymb,amsthm,amscd}
\chardef\bslash=`\\
\newtheorem{Thm}{Theorem}

\newtheorem{Lem}{Lemma}

\newtheorem{Rque}{Remark}

\newtheorem{Conj}{Conjecture}

\def\cit{{\mathbb C}}

\def\0{{\mathcal O}}
\def\g{{\mathfrak g}}
\def\p{{\mathfrak p}}

\begin{document}
\maketitle

\begin{abstract}
Based on our previous work \cite{Fu}, we prove that for any two projective symplectic resolutions $Z_1$ and $Z_2$ for a 
nilpotent orbit closure
 in a simple complex Lie algebra of classical type, then $Z_1$ is deformation equivalent to $Z_2$. This provides support for a 
``folklore'' conjecture on symplectic resolutions for symplectic singularities.
\end{abstract} 
\section{Introduction}
Let $X$ be a  complex variety, which is smooth in codimension 1. Following \cite{Bea}, the variety $X$  
is said to have {\em symplectic singularities} if there 
exists a holomorphic
symplectic 2-form $\omega$ on $X_{reg}$ such that for any resolution of singularities $\pi: Z \rightarrow X$, the 2-form
$\pi^* \omega$ defined {\em a priori} on $\pi^{-1}(X_{reg})$ can be extended to a holomorphic 2-form on $Z$. If furthermore 
the 2-form $\pi^* \omega$ extends to a holomorphic symplectic 2-form on the whole of $Z$ for some resolution of $X$, 
then we say that $X$ admits a {\em symplectic resolution}, and the resolution $\pi$ is called {\em symplectic}.

One important class of examples of symplectic singularities consists of nilpotent orbit closures in a semi-simple complex  Lie algebra. 
Projective symplectic resolutions for these singularities have been completely characterized in \cite{Fu}. This note studies the relationship
between two such projective symplectic resolutions.

Recall that two varieties $Z_1$ and $Z_2$ are {\em deformation equivalent} if there exists a flat morphism $\mathcal{Z} \xrightarrow{f} S$ over a 
connected curve $S$ such that $Z_1$ and $Z_2$ are isomorphic to two fibers of $f$. The purpose of this note is to prove the following:
\begin{Thm}\label{thm}
Let $\0$ be a nilpotent orbit in a complex simple Lie algebra $\g$ of classical type . Suppose that $\overline{\0}$ admits two projective
symplectic resolutions
$Z_1 \rightarrow \overline{\0}$ and $Z_2 \rightarrow \overline{\0}$. Then $Z_1$ is deformation equivalent to $Z_2$.  
\end{Thm}

In particular, we see that $Z_1$ and $Z_2$ are diffeomorphic, thus they have the same topological invariants (Betti numbers etc.), 
which is in some sense the McKay correspondence for nilpotent orbit closures.
Another motivation of this result is to provide support for the following:

\begin{Conj}
Let $X$ be an irreducible variety with symplectic singularities. Suppose that we have two projective symplectic resolutions:
$Z_1 \rightarrow X$ and $Z_2 \rightarrow X$. Then $Z_1$ is deformation equivalent to $Z_2$.
\end{Conj}

When $X$ is projective, this conjecture has been proven by D. Huybrechts (Theorem 4.6 \cite{Huy}). In \cite{Kal},
D. Kaledin proved the conjecture under a rather restrictive hypothesis (Condition 5.1 \cite{Kal}).
 In \cite{FN}, we  proved that any projective  symplectic 
resolution for a symmetric product $S^{(n)}$ of a symplectic connected 
surface $S$ is isomorphic to the Douady-Barlet resolution $S^{[n]} \rightarrow S^{(n)}$. 
In particular, this also verifies  the above conjecture. Some other results on uniqueness of symplectic resolutions for
some quotient symplectic singularities
and for some nilpotent orbit closures are obtained in \cite{FN}.

\section{Key lemma}
I am indebted to M. Brion for pointing out the following lemma, which plays a key role in the proof of our theorem.
\begin{Lem}\label{key}
Let $G$ be a semi-simple complex connected Lie group with Lie algebra $\g$. Let $P$ be a parabolic  subgroup of $G$ with Levi factor $L$.
Then $T^*(G/P)$ is deformation equivalent to $G/L$.
\end{Lem}
{\em Proof.}
 We first show that $T_o^*(G/P)$ is deformation equivalent to 
$P/L$ as $P$-varieties, where $o$ is a base point of $G/P$. Let $\p_u$ be the Lie algebra of the unipotent radical $P_u$ of $P$.
Then $\p_u$ equipped with the adjoint action of $P$ is identified with $T_o^*(G/P)$. 

Let $z$ be an element in the center of $\l:= \mathrm{Lie}(L)$ such that its centralizer in $G$ is exactly $L$, then its centralizer
in $\g$ is exactly $\l$. Consider the family $V: = (tz + \p_u)_{t\in \cit} \rightarrow{f} \cit$
 of sub-spaces in $\g$. Note that each subspace $tz + p_u$
is stable under the adjoint action of $L$ and  the adjoint action of $P_u$ (since $[z, \p_u]$ is contained in $\p_u$).
Thus we get a family of $P$-varieties with $\p_u$ being the special fiber.

Now we show that $tz + \p_u$ is isomorphic to $P/L$ for $t \neq 0$. In fact, the $P$-orbit of $tz$ coincides with the $P_u$-orbit of 
$t z$, which is closed in $\g$ (since $P_u$ is unipotent), but $Ad(P_u) tz$ and $tz + \p_u$ have the same tangent spaces at
$tz$, thus  $Ad(P_u) tz$ is also open in $tz + \p_u$, which shows that $tz + \p_u$ is exactly the $P$-orbit of $tz$. By our choice of $z$, the 
latter is isomorphic to $P/L$.

Now consider the family $G \times_P V \to \cit $, which is given by $(g,v)P \mapsto f(v)$. 
Then the central fiber is $T^*(G/P)$ and other fibers are all isomorphic to $G/L$, which concludes the proof.
 $\spadesuit$
\begin{Rque}\label{Rque1}
Notice that if two parabolic subgroups $P_1$ and $P_2$ have a Levi factor $L$ in commun, then the two deformations 
$G \times_{P_1} V_1 \to \cit$ and $G \times_{P_2} V_2 \to \cit$ have the same fibers over $\cit -\{0\}$. 
\end{Rque}
\begin{Rque}\label{Rque2}
The natural morphism $$G \times_P V \xrightarrow{\tilde{\pi}} W_P \subset \g, \quad (g, v) \mapsto Ad(g) v$$ gives a deformation 
of the Springer resolution $\pi: T^*(G/P) \to \overline{\0},$ where $W_P$ is the image of $\tilde{\pi}$, depending {\em a priori} on the
polarization $P$.
 Notice that $\tilde{\pi}_t$ is an isomorphism if $t \neq 0$.  
\end{Rque}

\section{Proof of the theorem}
First let us recall the following theorem from \cite{Fu}:
\begin{Thm}
Let $\0$ be a nilpotent orbit in a semi-simple complex Lie algebra $\g$ and $G$ a connected Lie group with Lie algebra $\g$. 
Suppose that $\overline{\0}$ admits a 
symplectic resolution $Z \rightarrow \overline{\0}$. Then $Z$ is isomorphic to $T^*(G/P)$ for some parabolic 
subgroup $P$ of $G$ and under this isomorphism, the map $Z \simeq T^*(G/P) \rightarrow \overline{\0}$ becomes the natural collapsing of the zero section.
\end{Thm}

Such a parabolic sub-group $P$ is called a {\em polarization} of $\0$ in \cite{Hes}.  
So to prove Theorem \ref{thm}, we need to show that if we have two polarizations $P_1, P_2$ of $\0$ such that $T^*(G/P_1)$
 and $T^*(G/P_2)$ are birational to $\0$, then $T^*(G/P_1)$ is deformation equivalent  to 
$T^*(G/P_2)$. In fact, we will prove that either $P_1$ and $P_2$ have conjugate  Levi factors or $G/P_1$ is isomorphic to $G/P_2$ (in some cases
of $\g = \mathfrak{so}_{2n}$). Then Lemma \ref{key} will conclude the proof.

To this end, we will do a case-by-case check, using the results of W. Hesselink in \cite{Hes}.
Let {\bf d}=$[d_1, \cdots, d_N] $ be the partition corresponding to the orbit $\0$ (c.f. section 5.1 \cite{CM}).
Let {\bf s} = $[s_1, \cdots, s_k]$ be the dual partition of {\bf d}, where $s_i = \#\{j| d_j \geq i \} $.
\subsection{Case $\g = \mathfrak{sl}_n$ }
Let $V = \cit^n$. A {\em flag} $F$ of $V$ is a sequence of sub-spaces $0 = F_0 \subsetneq F_1 \subsetneq \cdots \subsetneq F_k = V$. 
Its {\em type} is the sequence $(p_1, \cdots, p_k)$ given by $p_i = dim(F_i/F_{i-1})$. One shows easily that the stabilizer $P \subset G $ 
of $F$ is a parabolic subgroup. 

By Theorem 3.3 \cite{Hes}, any polarization $P$ of $\0$ is  a stabilizer of  such a flag $F$ with type $(s_{\sigma(1)}, \cdots, s_{\sigma(k)})$
for some permutation  $\sigma \in \mathfrak{S}_k$, and they have conjugate Levi factors.

\subsection{Case $\g = \mathfrak{sp}_{2n}$ }
Let $V = \cit^{2n}$ and $\phi$ a non-degenerate anti-symmetric bilinear form on $V$.
A flag $F = (F_0,\cdots,F_k)$ is called {\em isotropic} if $F_i^\bot = F_{k-i}$ for $0 \leq i \leq k$, where $F_i^\bot$ is the orthogonal
space of $F_i$ with respect to the bilinear form $\phi$. The type $(p_1, \cdots, p_k)$ of an isotropic flag satisfies $p_i = p_{k+1-i}$ for 
$1 \leq i \leq k$. By Lemma 4.3 \cite{Hes}, every parabolic subgroup $P$ of $G$ is the stabilizer of some isotropic flag, and
two parabolic subgroups with the same flag type are conjugate under the action of $G$ (Lemma 4.4 {\em loc. cit.}).

Here all congruences are modulo 2. 
For an even number $q$, let $$Pai(2n,q) = \{partitions \,\pi \,of\, 2n | \, \pi_j \equiv 1 \, if \, j \leq q;\, \pi_j \equiv 0 \, if \, j > q \}. $$
The union $Pai(2n) = \cup_q Pai(2n,q)$ parametrizes all conjugate classes of Levi factors of parabolic subgroups in $G$. For any $q$, there exists
an injective  Spaltenstein mapping (Proposition 6.5 \cite{Hes})
 $$ S_q: Pai(2n,q) \rightarrow \mathcal{P}_1(2n),  $$ where $\mathcal{P}_1(2n)$ is the set of partitions of $2n$ in 
which odd parts occur with even multiplicity (cf. Section 5.1 \cite{CM}). 

By the proof of Proposition 3.21 \cite{Fu}, if $P$ is a parabolic subgroup of $G$
such that $T^*(G/P)$ gives a symplectic resolution for $\overline{\0}$, then one has $q = \# \{j | d_j \equiv 1 \}$.
The injectivity of the map $S_q$ implies  that  any two such parabolic subgroups $P_1$ and $P_2$ have conjugate Levi factors.
\subsection{Case $\g = \mathfrak{so}_{2n+1}$}

Let $V = \cit^{2n+1}$ and let  $\phi$ be a non-degenerate symmetric bilinear form on $V$. In the same way as in the case of $\g = \mathfrak{sp}_{2n}$, one
defines the notion of isotropic flags. The proof goes in a similar way as in the case of $\g = \mathfrak{sp}_{2n}$.

\subsection{Case $\g = \mathfrak{so}_{2n}$}

Let $V = \cit^{2n}$ and $\phi$ a non-degenerate symmetric bilinear form on $V$. The group $H$ of  automorphisms of $V$ preserving $\phi$
has two components. The identity component of $H$ is our Lie group $G \simeq SO(2n)$. 

By Lemma 4.4 \cite{Hes}, the class of parabolic subgroups of $G$ 
with flag type $(p_1, \cdots, p_{k})$ splits into two conjugacy classes (denoted by $P_1$ and $P_2$) under the action of $G$ if and only if  
$k = 2 t$ and  $p_t \geq 2$. Furthermore, the two parabolic subgroups are conjugate under the action of $H$, i.e. there exists an element 
$h \in H$ such that $P_2 = hP_1h^{-1}$. Take an isotropic flag $F = (F_0, \cdots, F_k) \in G/P_1$, 
then $F' = hF = (hF_0, \cdots, hF_k)$ is an isotropic flag in $G/P_2$. This gives an isomorphism 
between $G/P_1$ and $G/P_2$. 

So we need only to consider polarizations of $\0$ with different flag types. By the proof of Proposition 3.22 \cite{Fu} and the proof of
Lemma 4.6 \cite{Hes}, two such polarizations have conjugate Levi factors.. 

\section{Concluding remarks} 

 {\bf 4.1} Though we believe that Theorem \ref{thm} is also true for exceptional Lie algebras, by lack of a complete description of 
polarizations of their nilpotent orbits, we don't know how to check this.\\

{\bf 4.2}  One should bear in mind that for two symplectic resolutions $Z_i \to \overline{\0}, i=1,2,$ though $Z_1$ is deformation
equivalent to $Z_2$, $Z_1$ and $Z_2$ may be non-isomorphic. An explicit example is given in \cite{FN}. \\

{\bf 4.3} By Remark \ref{Rque1}, our theorem can be strengthened as follows: Let $Z_i \to \overline{\0}, i=1,2,$ be two symplectic resolutions
for a nilpotent orbit closure in a simple complex classical Lie algebra. Then there exists two deformations $\mathcal{Z}_i \to \cit,$
such that $\mathcal{Z}_{i,0} \cong Z_i$ and $\mathcal{Z}_{1,t} \cong \mathcal{Z}_{2,t} $ for any $t \neq 0$. This can be viewed as an analogue
of the result of D. Huybrechts (Theorem 4.6 \cite{Huy}). \\

{\bf 4.4}  A stronger form of Conjecture 1 is stated in \cite{FN}, where we conjectured that there
exist deformations $F_i: \mathcal{Z}_i \to \mathcal{X}$   of the morphisms $\pi_i: Z_i \to X$, such that $F_{i,t}$ is an isomorphism
for $t \neq 0$. 
In the case of nilpotent orbit closures in $\mathfrak{sl}_n$, we proved this in \cite{FN},  
by constructing  explicitly the deformations. For the other cases, though by Remark 2 we have a deformation $\tilde{\pi}$ of $\pi$, 
we don't know whether the deformation spaces $W_P$ of $\overline{\0}$ are isomorphic or not for two polarizations.  \\

{\em Acknowledgments:} I want to thank A. Beauville, M. Brion and Y. Namikawa for helpful discussions and suggestions. 
It was pointed out to me by Y. Namikawa that there was an error in an earlier version of this note. Then M. Brion suggested
to me the Key lemma to complete the proof.

\quad \\
Baohua Fu,\\
Labortoire J. A. Dieudonn\'e, Parc Valrose \\ 06108 Nice cedex 02, FRANCE \\
baohua.fu@polytechnique.org 

\end{document}